\documentclass[12pt,reqno]{amsart} 

\usepackage[dvips]{graphicx} 
\usepackage[arrow,matrix,curve]{xy}

\usepackage{amssymb,latexsym, amsmath, amscd, array, hyperref}

\usepackage{tikz}     

                            {\egroup\par\bigskip}

\def\tag#1#2{\hbox to\textwidth{#1\hfil$\displaystyle #2$\hfil}}

\newtheorem{theorem}{Theorem}[section]

\newtheorem{corollary}[theorem]{Corollary}

\theoremstyle{definition}
\newtheorem{definition}[theorem]{Definition}
\newtheorem{remark}[theorem]{Remark}

\newcommand{\ww}{{\sqrt{\tfrac{g}{\ell}}}}

\DeclareMathOperator{\adequal}{\;\raisebox{-3pt}{$\ulcorner\!\urcorner$}\;}

\DeclareMathOperator{\N}{{\mathbb N}}

\DeclareMathOperator{\Q}{{\mathbb Q}}
\DeclareMathOperator{\R}{{\mathbb R}}
\DeclareMathOperator{\C}{{\mathbb C}}

\newcommand\astr{{{}^\ast\!\R}} 

\newcommand\astn{{{}^\ast\!\N}} 

\newcommand\astc{{{}^{\ast}\!\C}}

\numberwithin{figure}{section} \numberwithin{table}{section}

\begin{document}

\title [Small oscillations of the pendulum] {Small oscillations of the
pendulum, Euler's method, and adequality}

\author{Vladimir Kanovei} \address{V. Kanovei, IPPI, Moscow, and MIIT,
Moscow, Russia}\email{kanovei@googlemail.com}

\author[K. Katz]{Karin U. Katz}\address{K. Katz, Department of
Mathematics, Bar Ilan University, Ramat Gan 52900 Israel}
\email{katzmik@math.biu.ac.il}

\author[Mikhail Katz]{Mikhail G. Katz} \address{M. Katz, Department of
Mathematics, Bar Ilan University, Ramat Gan 52900 Israel}
\email{katzmik@macs.biu.ac.il}

\author{Tahl Nowik}\address{T. Nowik, Department of Mathematics, Bar
Ilan University, Ramat Gan 52900 Israel}\email{tahl@math.biu.ac.il}

\begin{abstract}
Small oscillations evolved a great deal from Klein to Robinson.  We
propose a concept of solution of differential equation based on
Euler's method with infinitesimal mesh, with well-posedness based on a
relation of adequality following Fermat and Leibniz.  The result is
that the period of infinitesimal oscillations is independent of their
amplitude.

\end{abstract}

\maketitle

\tableofcontents

\section{Small oscillations of a pendulum}

The breakdown of infinite divisibility at quantum scales makes
irrelevant the mathematical definitions of derivatives and integrals
in terms of limits as~$x$ tends to zero.  Rather, quotients like
$\frac{\Delta y}{\Delta x}$ need to be taken in a certain range, or
level.  The work \cite{NK} developed a general framework for
differential geometry at level~$\lambda$, where~$\lambda$ is an
infinitesimal but the formalism is a better match for a situation
where infinite divisibility fails and a scale for calculations needs
to be fixed accordingly.  In this paper we implement such an approach
to give a first rigorous account ``at level~$\lambda$'' for small
oscillations of the pendulum.

In his 1908 book \emph{Elementary Mathematics from an Advanced
Standpoint}, Felix Klein advocated the introduction of calculus into
the highschool curriculum.  One of his arguments was based on the
problem of small oscillations of the pendulum.  The problem had been
treated until then using a somewhat mysterious \emph{superposition
principle} involving a hypothetical circular motion of the pendulum.
Klein advocated what he felt was a better approach, involving the
differential equation of the pendulum; see \cite[p.\;187]{Kl08}.

The classical problem of the pendulum translates into the second order
nonlinear differential equation~$\ddot x=-\frac{g}{\ell}\sin x$ for
the variable angle~$x$ with the vertical direction, where~$g$ is the
constant of gravity and~$\ell$ is the length of the (massless) rod or
string.  The problem of small oscillations deals with the case of
small amplitude, i.e.,~$x$ is small, so that~$\sin x$ is
approximately~$x$.

Then the equation is boldly replaced by the linear one~$\ddot
x=-\frac{g}{\ell}x$, whose solution is harmonic motion with
period~$2\pi\sqrt{\ell/g}$.

This suggests that the period of small oscillations should be
independent of their amplitude.  The intuitive solution outlined above
may be acceptable to a physicist, or at least to the mathematicians'
proverbial physicist.  The solution Klein outlined in his book does
not go beyond the physicist's solution.

The Hartman--Grobman theorem \cite{Ha60}, \cite{Gr59} provides a
criterion for the flow of the nonlinear system to be conjugate to that
of the linearized system, under the hypothesis that the linearized
matrix has no eigenvalue with vanishing real part.  However, the
hypothesis is not satisfied for the pendulum problem.

To give a rigorous mathematical treatment, it is tempting to exploit a
hyperreal framework following \cite{NK}.  Here the notion of small
oscillation can be given a precise sense, namely oscillation with
infinitesimal amplitude.  

However even for infinitesimal~$x$ one cannot boldly replace~$\sin x$
by~$x$.  Therefore additional arguments are required.

The linearisation of the pendulum is treated in \cite{St15} using
Dieners' ``Short Shadow'' Theorem; see Theorem~5.3.3 and Example~5.3.4
there.  This text can be viewed as a self-contained treatement of
Stroyan's Example~5.3.4.

The traditional setting exploiting the real continuum is only able to
make sense of the claim that \emph{the period of small oscillations is
independent of the amplitude} by means of a paraphrase in terms of
limits.  In the context of an infinitesimal-enriched continuum, such a
claim can be formalized more literally; see Corollary~\ref{c1061}.
What enables us to make such distinctions is the richer syntax
available in Robinson's framework.

Terence Tao has recently authored a number of works exploiting
ultraproducts in general, and Robinson's infinitesimals in particular,
as a fundamental tool; see e.g., \cite{Ta14}, \cite{TV}.  In the
present text, we apply such an approach to small oscillations.

Related techniques were exploited in \cite{LS}.

\section
{Vector fields, walks, and integral curves}
\label{s2}

The framework developed in \cite{Ro66} involves a proper extension
$\astr\supseteq\R$ preserving the properties of~$\R$ to a large extent
discussed in Remark~\ref{r22}.  Elements of~$\astr$ are called
hyperreal numbers.  A positive hyperreal number is called
\emph{infinitesimal} if it is smaller than every positive real number.

We choose a fixed positive infinitesimal~$\lambda\in\astr$ (a
restriction on the choice of~$\lambda$ appears in Section~\ref{s7}).
Given a classical vector field~$V=V(z)$ where~$z\in\C$, one forms an
infinitesimal displacement~$\delta_ F(z)=\lambda V(z)$ with the aim of
constructing the integral curves of the corresponding flow~$F_t$ in
the plane.  Note that a zero of~$\delta_ F$ corresponds to a fixed
point (i.e., a constant integral ``curve'') of the flow.  The
infinitesimal generator is the function~$F:\astc\to\astc$, also called
a \emph{prevector field}, defined by
\begin{equation}
\label{e21}
F(z)=z+\delta_ F(z),
\end{equation}
where~$\delta_ F(z)=\lambda V(z)$ in the case of a displacement
generated by a classical vector field as above, but could be a more
general internal function~$F$ as discussed in \cite{NK}.

We propose a concept of solution of differential equation based on
Euler's method with infinitesimal step size, with well-posedness based
on a property of adequality (see Section~\ref{s3}), as follows.

\begin{definition}
The \emph{hyperreal flow}, or \emph{walk},~$F_t(z)$ is a
$t$-parametrized map~$\astc\to\astc$ defined whenever~$t$ is a
hypernatural multiple~$t=N\lambda$ of~$\lambda$, by setting
\begin{equation}
\label{e22}
F_t(z)=F_{N\lambda}^{\phantom{I}}(z)=F^{\circ N}(z),
\end{equation}
where~$F^{\circ N}$ is the~$N$-fold composition.
\end{definition}

The fact that the infinitesimal generator~$F$ given by~\eqref{e21} is
invariant under the flow~$F_t$ of~\eqref{e22} receives a transparent
meaning in this framework, expressed by the commutation relation
$F\circ F^{\circ N}=F^{\circ N} \circ F$ due to \emph{transfer} (see
Remark~\ref{r22}) of associativity of composition of maps.

\begin{remark}
\label{r22}
The \emph{transfer principle} is a type of theorem that, depending on
the context, asserts that rules, laws or procedures valid for a
certain number system, still apply (i.e., are \emph{transfered}) to an
extended number system.  Thus, the familiar extension~$\Q\subseteq\R$
preserves the property of being an ordered field.  To give a negative
example, the extension~$\R\subseteq\R\cup\{\pm\infty\}$ of the real
numbers to the so-called \emph{extended reals} does not preserve the
property of being an ordered field.  The hyperreal extension
$\R\subseteq\astr$ preserves \emph{all} first-order properties,
including the identity~$\sin^2 x + \cos^2 x =1$ (valid for all
hyperreal~$x$, including infinitesimal and infinite values
of~$x\in\astr$).  The natural numbers~$\N\subseteq\R$ are naturally
extended to the hypernaturals~$\astn\subseteq\astr$.  For a more
detailed discussion, see the textbook \emph{Elementary Calculus}
\cite{Ke86}.
\end{remark}

\begin{definition}
The \emph{real flow}~$f_t$ on~$\C$ for~$t\in\R$ when it exists is
constructed as the shadow (i.e., standard part) of the hyperreal
walk~$F_t$ by setting
$f_t(z)=\text{st}\left(F_{N\lambda}^{\phantom{I}}(z)\right)$ where
$N=\left\lfloor \frac{t}{\lambda}\right\rfloor$, while~$\lfloor
x\rfloor$ rounds off the number~$x$ to the nearest hyperinteger no
greater than~$x$, and ``st'' (standard part or shadow) rounds off each
finite hyperreal to its nearest real number.
\end{definition}

For~$t$ sufficiently small, suitable regularity conditions ensure that
the point~$F_{N\lambda}(z)$ is finite so that the shadow is defined.

The usual relation of being infinitely close is denoted~$\approx$.
Thus~$z,w$ satisfy~$z\approx w$ if and only if~$\text{st}(z-w)=0$.
This relation is an additive one (i.e., invariant under addition of a
constant).

The appropriate relation for working with small prevector fields is
not additive but rather multiplicative (i.e., invariant under
multiplication by a constant), as detailed in Section~\ref{s3}.

\section
{Adequality}
\label{s3}

We will use Leibniz's notation~$\adequal$ to denote the relation of
adequality (see below).  Leibniz actually used a symbol that looks
more like~$\sqcap$ but the latter is commonly used to denote a
product.  Leibniz used the symbol to denote a generalized notion of
equality ``up to'' (though he did not distinguish it from the usual
symbol~$=$ which he also used in the same sense).  A prototype of such
a relation (though not the notation) appeared already in Fermat under
the name \emph{adequality}.  For a re-appraisal of Fermat's
contribution to the calculus see \cite{KSS}; for Leibniz's, see
\cite{KS1}, \cite{Ba16b}; for Euler see \cite{B11}, \cite{Ba16}; for
Cauchy's contribution, see \cite{KK11b}, \cite{KK12a}, \cite{BK},
\cite{Ba14}.  We will use the sign~$\adequal$ for a multiplicatively
invariant relation among (pre)vectors defined as follows.

\begin{definition}
Let~$z,w\in\astc$.  We say that~$z$ and~$w$ are \emph{adequal} and
write~$ z\adequal w~$ if either~$\frac{z}{w}\approx 1$ (i.e.,
$\frac{z}{w}- 1$ is infinitesimal) or~$z=w=0$.
\end{definition}
This implies in particular that the angle between~$z,w$ (when they are
nonzero) is infinitesimal, but~$\adequal$ is a stronger condition.  If
one of the numbers is appreciable, then so is the other and the
relation~$z\adequal w$ is equivalent to~$z\approx w$.  If one of~$z,w$
is infinitesimal then so is the other, and the difference~$|z-w|$ is
not merely infinitesimal, but so small that the quotients~$|z-w|/z$
and~$|z-w|/w$ are infinitesimal, as well.

We are interested in the behavior of orbits in a neighborhood of a
fixed point~$0$, under the assumption that the infinitesimal
displacement satisfies the Lipschitz condition.  In such a situation,
we have the following theorem.

\begin{theorem}
\label{t25}
Assume that for some finite~$K$, we
have~$\delta_F(z)-\delta_F(w)<K\lambda|z-w|$.  Then prevector fields
defined by adequal infinitesimal displacements produce hyperreal walks
that are adequal at each finite time, or in formulas: if~$\delta_ F
\adequal \delta_ G$ then~$F_t \adequal G_t$ when~$t$ is finite.
\end{theorem}

This was shown in \cite[Example~5.12]{NK}.

\section{Infinitesimal oscillations}

Let~$x$ denote the variable angle between an oscillating pendulum and
the downward vertical direction.  By considering the projection of the
force of gravity in the direction of motion, one obtains the equation
of motion~$ m\ell \ddot x = -mg \sin x~$ where~$m$ is the mass of the
bob of the pendulum,~$\ell$ is the length of its massless rod, and~$g$
is the constant of gravity.  Thus we have a second order nonlinear
differential equation
\begin{equation}
\label{e1061b}
\ddot x=- \tfrac{g}{\ell} \sin x.
\end{equation}
The initial condition of releasing the pendulum at angle~$a$ (for
\emph{amplitude}) is
\[
\begin{cases}
x(0) = a, \\ \dot x(0)=0.
\end{cases}
\]
We replace \eqref{e1061b} by the pair of
first order equations
\[
\begin{cases}
\dot x = \ww\, y, \\ \dot y =- \ww \sin x,
\end{cases}
\]
and initial condition~$(x,y)=(a,0)$.  We identify~$(x,y)$
with~$z=x+iy$ and~$(a,0)$ with~$a+i0$ as in Section~\ref{s2}.  The
classical vector field corresponding to this system is then
\begin{equation}
\label{e1061}
X(x, y) = \ww\, y - i \ww \sin x.
\end{equation}
The corresponding prevector field~$F$ is defined by the infinitesimal
displacement~$\delta_ F(z) = \lambda \ww \, y - i\lambda \ww \sin x$
so that~$F(z)=z+\delta_ F(z)$.  We are interested in the flow of~$F$,
with initial condition~$a+ 0i$, generated by hyperfinite iteration
of~$F$.

Consider also the linearisation, i.e., prevector field~$E(z)=z+\delta_
E(z)$ defined by the displacement
\[
\delta_ E(z) = \lambda\ww\, y - i\lambda\ww\, x= -i \lambda\ww z
\]
where as before~$z=x+iy$.  We are interested in small oscillations,
i.e., the case of infinitesimal amplitude~$a$.  Since sine is
asymptotic to the identity function for infinitesimal inputs, we
have~$\delta_ E\adequal \delta_ F.$ Due to the multiplicative nature
of this relation, the rescalings of~$E$ and~$F$ by change of
variable~$z=aZ$ remain adequal and therefore define adequal walks and
identical real flows by Theorem~\ref{t25}.

\section{Adjusting linear prevector field}

We will compare~$E$ to another linear prevector field
\[
\begin{aligned}
H(x&+iy)=\\ & =e^{-i\lambda\ww} (x+iy) \\& =
\left(x\cos\lambda\ww+y\sin\lambda\ww\,\right)+
\left(-x\sin\lambda\ww+y\cos\lambda\ww\,\right)i
\end{aligned}
\]
given by clockwise rotation of the~$x,y$ plane by infinitesimal
angle~$\lambda\ww$, so that
\[
\delta_H^{\phantom{I}}(z)=\left(e^{-i\lambda\ww}-1\right)z.
\]

The corresponding hyperreal walk, defined by hyperfinite iteration
of~$H$, satisfies the exact equality
\begin{equation}
\label{e1062b}
H_t(a, 0)= \left(a \cos \ww\, t, - a \sin \ww\, t\right)
\end{equation}
whenever~$t$ is a hypernatural multiple of~$\lambda$.  In particular,
we have the periodicity property~$H_{\frac{2\pi}{\sqrt{g/\ell}}}(z)=z$
and therefore
\begin{equation}
\label{e1064}
H_{t+\frac{2\pi}{\sqrt{g/\ell}}}=H_t
\end{equation}
whenever both~$t$ and~$\frac{2\pi}{\sqrt{g/\ell}}$ are hypernatural
multiples of~$\lambda$.  Note that we have~$\delta_E(z)=-i\lambda\ww
z$ and $\delta_H(z)=\big(e^{-i\lambda\ww}-1\big)z$.  Therefore
\[
\frac{\delta_H(z)}{\delta_E(z)} =
\frac{e^{-i\lambda\ww}-1}{-i\lambda\ww}=
\frac{1}{\lambda\ww}\Big(\sin\lambda\ww + \big(\cos \lambda\ww
-1\big)i \Big)
\]
The usual estimates give
\begin{equation}
\label{e1062}
\frac{\sin\lambda\ww}{\lambda\ww}\approx 1, \quad \frac{\cos
\lambda\ww -1}{\lambda\ww}\approx 0
\end{equation}
so~$\frac{\delta_H(z)}{\delta_E(z)} \approx 1 +0i=1$, that is
$\delta_H(z) \adequal \delta_E(z)$.  By Theorem~\ref{t25}
%
%
the hyperfinite walks of~$E$ and~$H$
satisfy~$E_t(a,0)\adequal{}H_t(a,0)$ for each finite initial
amplitude~$a$ and for all finite time~$t$ which is a hypernatural
multiple of~$\lambda$.

\section{Conclusion}
\label{s7}

The advantage of the prevector field~$H$ is that its hyperreal walk is
given by an explicit formula \eqref{e1062b} and is therefore periodic
with period precisely~$\frac{2\pi}{\ww}$, provided we choose our base
infinitesimal~$\lambda$ in such a way that~$\frac{2\pi}{\lambda \ww}$
is hypernatural.

We obtain the following consequence of~\eqref{e1064}: modulo a
suitable choice of a representing prevector field (namely,~$H$) in the
adequality class, the hyperreal walk is periodic with period
$2\pi\sqrt{\ell/g}$.  This can be summarized as follows.

\begin{corollary}
\label{c1061}
The period of infinitesimal oscillations of the pendulum is
independent of their amplitude.  
\end{corollary}

If one rescales such an infinitesimal oscillation to appreciable size
by a change of variable~$z= a Z$ where~$a$ is the amplitude, and takes
standard part, one obtains a standard harmonic oscillation with
period~$2\pi\sqrt{\ell/g}$.  The formulation contained in
Corollary~\ref{c1061} has the advantage of involving neither rescaling
nor shadow-taking.

\section*{Acknowledgments} 
We are grateful to Jeremy Schiff for drawing our attention to the
Hartman--Grobman theorem, and to Semen Kutateladze and Dalibor Pra\v
z\'ak for some helpful suggestions.  M. Katz was partially supported
by the Israel Science Foundation grant no.~1517/12.

\bibliographystyle{amsalpha}

\begin{thebibliography}{AI}

\bibitem[Bair et al.~2013]{B11} Bair, J.; B\l{}aszczyk, P.; Ely, R.;
Henry, V.; Kanovei, V.; Katz, K.; Katz, M.; Kutateladze, S.; McGaffey,
T.; Schaps, D.; Sherry, D.; Shnider,~S.  ``Is mathematical history
written by the victors?''  \emph{Notices of the American Mathematical
Society} \textbf{60}, no.~7, 886-904.

See \url{http://www.ams.org/notices/201307/rnoti-p886.pdf}

and \url{http://arxiv.org/abs/1306.5973}




\bibitem[Bair et al.~2016]{Ba16} Bair, J.; B\l{}aszczyk, P.; Ely, R.;
Henry, V.; Kanovei, V.; Katz, K.; Katz, M.; Kutateladze, S.; McGaffey,
T.; Reeder, P.; Schaps, D.; Sherry, D.; Shnider, S.  ``Interpreting
the infinitesimal mathematics of Leibniz and Euler.''  \emph{Journal
for General Philosophy of Science} (2016), to appear.


\bibitem[Bascelli et al. 2014]{Ba14} Bascelli, T.; Bottazzi, E.;
Herzberg, F.; Kanovei, V.; Katz, K.; Katz, M.; Nowik, T.; Sherry, D.;
Shnider, S.  ``Fermat, Leibniz, Euler, and the gang: The true history
of the concepts of limit and shadow.''  \emph{Notices of the American
Mathematical Society} \textbf{61}, no.~8, 848-864.


\bibitem[Bascelli et al.~2016]{Ba16b} Bascelli, T.; B\l{}aszczyk, P.;
Kanovei, V.; Katz, K.; Katz, M.; Schaps, D.; Sherry, D.  ``Leibniz vs
Ishiguro: Closing a quarter-century of syncategoremania.''
\emph{HOPOS: Journal of the Internatonal Society for the History of
Philosophy of Science} \textbf{6}, no.~1.  See
\url{http://dx.doi.org/10.1086/685645}

and \url{http://arxiv.org/abs/1603.07209}







\bibitem[Borovik \& Katz 2012]{BK} Borovik, A., Katz, M.  ``Who gave
you the Cauchy--Weier\-strass tale?  The dual history of rigorous
calculus.''  \emph{Foundations of Science} \textbf{17}, no.~3,
245-276.  See \url{http://dx.doi.org/10.1007/s10699-011-9235-x}





\bibitem[Grobman 1959]{Gr59} Grobman, D.  ``Homeomorphisms of systems
of differential equations.''  \emph{Doklady Akademii Nauk SSSR}
\textbf{128}, 880--881.

\bibitem[Hartman 1960]{Ha60} Hartman, P.  ``A lemma in the theory of
structural stability of differential equations.''  \emph{Proceedings
of the American Mathematical Society} \textbf{11} (4), 610--620.


\bibitem[Katz \& Katz 2011]{KK11b} Katz, K.; Katz, M.  ``Cauchy's
continuum.''  \emph{Perspectives on Science} \textbf{19}, no.~4,
426-452.  See \url{http://arxiv.org/abs/1108.4201} and

\url{http://www.mitpressjournals.org/doi/abs/10.1162/POSC_a_00047}


\bibitem[Katz \& Katz 2012]{KK12a} Katz, K.; Katz, M.  ``A Burgessian
critique of nominalistic tendencies in contemporary mathematics and
its historiography.''  \emph{Foundations of Science} \textbf{17},
no.~1, 51--89.  

See \url{http://dx.doi.org/10.1007/s10699-011-9223-1}

and \url{http://arxiv.org/abs/1104.0375}



\bibitem[Katz, Schaps \& Shnider 2013]{KSS} Katz, M.; Schaps, D.;
Shnider, S.  ``Almost Equal: The Method of Adequality from Diophantus
to Fermat and Beyond.''  \emph{Perspectives on Science} \textbf{21},
no.~3, 283-324.


\bibitem[Katz \& Sherry 2013]{KS1} Katz, M., Sherry, D.  ``Leibniz's
infinitesimals: Their fictionality, their modern implementations, and
their foes from Berkeley to Russell and beyond.''  \emph{Erkenntnis}
\textbf{78}, no.~3, 571--625.  

See \url{http://dx.doi.org/10.1007/s10670-012-9370-y}

and \url{http://arxiv.org/abs/1205.0174}





\bibitem[Keisler 1986]{Ke86} Keisler, H. J.  \emph{Elementary
Calculus: An Infinitesimal Approach.}  Second Edition.  Prindle, Weber
\& Schimidt, Boston.  See online version at
\url{http://www.math.wisc.edu/~keisler/calc.html}


\bibitem[Klein 1908]{Kl08} Klein, F.  \emph{Elementary Mathematics
from an Advanced Standpoint.  Vol.~I.  Arithmetic, Algebra, Analysis.}
Translation by E. R. Hedrick and C. A. Noble [Macmillan, New York,
1932] from the third German edition [Springer, Berlin, 1924].
Originally published as \emph{Elementarmathematik vom h\"oheren
Standpunkte aus} (Leipzig, 1908).

\bibitem[Lobry \& Sari 2008]{LS} Lobry, C.; Sari, T.  ``Non-standard
analysis and representation of reality.''  \emph{Internat. J. Control}
\textbf{81}, no.~3, 517--534.


\bibitem[Nowik \& Katz 2015]{NK} Nowik, T., Katz, M.  ``Differential
geometry via infinitesimal displacements.''  \emph{Journal of Logic
and Analysis} \textbf{7}:5, 1--44.  See

\url{http://www.logicandanalysis.org/index.php/jla/article/view/237/106}
and \url{http://arxiv.org/abs/1405.0984}



\bibitem[Pra\v z\'ak et al.~2016]{Pr16} Pra\v z\'ak, D.; Rajagopal,
K.; Slav\'\i k, J.  ``A non-standard approach to a constrained forced
oscillator.''  Preprint.
%
%

\bibitem[Robinson 1966]{Ro66} Robinson, A.  \emph{Non-standard
analysis}.  North-Holland Publishing, Amsterdam.

\bibitem[Stroyan 2015]{St15} Stroyan, K.  \emph{Advanced Calculus
using Mathematica: NoteBook Edition}.

\bibitem[Tao 2014]{Ta14} Tao, T.  \emph{Hilbert's fifth problem and
related topics}.  Graduate Studies in Mathematics, 153.  American
Mathematical Society, Providence.

\bibitem[Tao \& Vu 2016]{TV} Tao, T.; Van Vu, V.  ``Sum-avoiding sets
in groups.''  

See \url{http://arxiv.org/abs/1603.03068}


\end{thebibliography}

\end{document}